\documentclass[12pt]{amsart}

\usepackage{amssymb}
\usepackage{bm}
\usepackage{graphicx}
\usepackage{psfrag}
\usepackage[usenames,dvipsnames]{xcolor}
\usepackage{subcaption}
\usepackage{enumerate}
\usepackage{url}

\usepackage{algpseudocode}

\usepackage{mathdots}
\usepackage{mathtools}

\usepackage{xy}
\input xy
\xyoption{all}

\usepackage{tikz}
\usetikzlibrary{positioning} 
\usetikzlibrary{graphs,graphs.standard}

\newcommand{\excise}[1]{}

\newtheorem*{firstisothm}{The First Isomorphism Theorem}
\newtheorem*{freshmansdream}{The Freshman's Dream}
\newtheorem*{rootthm}{The Root Theorem}
\newtheorem*{keylemma}{The Key Lemma}
\newtheorem*{fundamentalthm}{The Fundamental Theorem of Finite Fields}

\theoremstyle{definition}

\numberwithin{equation}{section}

\newcommand{\ring}[1]{\ensuremath{\mathbb{#1}}}

\renewcommand\>{\rangle}

\newcommand\KK{{\mathbb K}}

\newcommand\ZZ{\ring{Z}}
\newcommand\FF{\ring{F}}

\newcommand\iso{\cong}

\def\ol#1{{\overline {#1\vphantom{1}}}}

\DeclareMathOperator\image{Im} 
\DeclareMathOperator\lcm{lcm} 

\begin{document}

\mbox{}
\title[Fundamental theorem of finite fields]{The fundamental theorem of finite fields: \\ a proof from first principles}

\author[A.~Chavez]{Anastasia Chavez}
\address{Mathematics Department\\University of California Davis\\Davis, CA 95616}
\email{amrchavez@ucdavis.edu}

\author[C.~O'Neill]{Christopher O'Neill}
\address{Mathematics and Statistics Department\\San Diego State University\\San Diego, CA 92182}
\email{cdoneill@sdsu.edu}

\date{\today}

\begin{abstract}
A mathematics student's first introduction to the fundamental theorem of finite fields (FTFF) often occurs in an advanced abstract algebra course and invokes the power of Galois theory to prove it. Yet the combinatorial and algebraic coding theory applications of finite fields can show up early on for students in STEM. To make the FTFF more accessible to students lacking exposure to Galois theory, we provide a proof from algebraic ``first principles.''
\end{abstract}

\maketitle


\section{Introduction}
\label{sec:intro}

A student's first introduction to finite fields and the magic they invoke often occurs in an advanced undergraduate or graduate abstract algebra course. In particular, the fundamental theorem of finite fields (FTFF) is most commonly proved via Galois theory. Finite fields have many exciting combinatorial applications, one of which is algebraic coding theory. Error-correcting codes, t-designs, and Hamming codes are common topics for computer science majors with minimal abstract algebra training. For curious undergraduates with just one year of abstract algebra, such applied combinatorics is both enriching and inspiring. Yet the Galois theory approach to finite fields leaves these students at a disadvantage. To fill an apparent gap in the accessibility of the FTFF, we provide a proof of this great theorem from ``first principles," i.e., without appealing to Galois groups or splitting fields.

\begin{fundamentalthm}\label{t:ftff}
\text{ }
\begin{enumerate}[(a)]
\item 
There is a field with exactly $q$ elements if and only if $q = p^r$ for $p$ prime,~$r \ge 1$.  

\item 
Any two finite fields of the same cardinality are isomorphic.  

\item 
For any finite field $\FF$ with $|\FF| = p^r$ for $p$ prime, 
\begin{enumerate}[(i)]
\item 
the additive group $(\FF, +) \cong ((\ZZ_p)^r, +)$, and

\item 
the multiplicative group $(\FF \setminus \{0\}, \cdot)$ is cyclic.  

\end{enumerate}
\end{enumerate}
\end{fundamentalthm}

The proof we provide here is built from several different sources, many of which either briefly mention or wave their hands at Galois theory for one or more parts of the argument \cite{biggs,dummitfoote,hungerford,lang}.  This approach of introducing finite fields to those with little abstract algebra exposure has been successfully implemented in several iterations of the applied combinatorics course at our former home institution.  
We provide this manuscript as a resource for those in need of a proof of the Fundamental Theorem that does not utilize the heavy machinery of Galois theory.  

This article is organized as follows. In Section~\ref{sec:background}, we survey the assumed abstract algebra background, and in Section~\ref{sec:constructing}, we review quotient rings and outline a general method for explicitly constructing finite fields.  Sections~\ref{sec:keylemma} and~\ref{sec:mainproof} together contain the proof of the FTFF, with the former section providing a Key Lemma that has some consequences of its own. 

%
%
%
%
%
%
%
%

\subsection*{Acknowledgements}

The authors would like to thank Scott Chapman, Lily Silverstein and Wencin Poh for numerous helpful conversations.  We are grateful to Jes\'{u}s~A. De~Loera for sharing his lecture notes from the course that inspired this manuscript.  We also wish to thank the referees for their insightful comments.

\section{Prerequisite background}
\label{sec:background}

In this section, we survey the minimal prerequisite definitions and results that are needed for this article.  We assume the reader is familiar with undergraduate-level linear algebra (including vector space dimension) and ring theory (including cosets and quotient rings). Please see \cite{dummitfoote,hungerford} for more details.  Note that all rings are assumed to be commutative and have a multiplicative identity.  

There are two main families of rings appearing in this paper.  The first is the ring~$\ZZ_n$ of integers modulo $n \ge 2$.  Note that $\ZZ_n$ is a field whenever $n$ is prime, and contains zero-divisors whenever $n$ is composite.  The second is the polynomial ring $F[x]$ whose coefficient ring $F$ is a field, as well as quotients $F[x]/I$ by an ideal $I$.  Several times throughout this article, we will use the fact that the polynomial ring $F[x]$ is:
\begin{enumerate}[(i)]
\item 
a \emph{principal ideal domain (PID)}, meaning every ideal $I \subset F[x]$ can be written as $I = \langle f(x) \rangle$ for some $f(x) \in I$, and the quotient ring $F[x]/I$ is a field if and only if $f(x)$ is irreducible; and

\item 
a \emph{unique factorization domain (UFD)}, meaning that every monic, non constant polynomial in $F[x]$ can be written uniquely (up to reordering) as a product of monic irreducible polynomials in $F[x]$.  

\end{enumerate}

We close this section with the following theorems, all of which are usually covered in an introductory course in rings, and which will be used in the proof of the Key Lemma given in Section~\ref{sec:keylemma}.   

\begin{rootthm}
Fix a field $F$.  For any $f \in F[x]$, we have $f(a) = 0$ if and only if $f(x) = (x - a)g(x)$ for some $g \in F[x]$.   
\end{rootthm}

\begin{freshmansdream}
Fix a prime $p$, and let $R$ be a ring with characteristic~$p$ (that is, $p \cdot 1 = 0$ in $R$).  If $a, b \in R$, then $(a + b)^p = a^p + b^p$.  
\end{freshmansdream}


\begin{firstisothm}
If $R$ and $S$ are rings and $\sigma: R \rightarrow S$ is a ring homo\-morphism, then $\image(\sigma) \subset S$ is a subring, $\ker(\sigma) \subset R$ is an ideal, and $R/\ker{\sigma} \cong \image(\sigma)$.
\end{firstisothm}

\section{Constructing finite fields}
\label{sec:constructing}

When constructing small finite fields from first principles, a common approach is to use the addition and multiplication tables (or ``$+/\cdot$ tables'') to help guide the behavior of the field's elements.  For example, suppose we wish to discover all possible finite fields with $3$ elements.  We know there must be two distinguished elements $0$ and $1$, so denoting the only remaining element $a$, we can consider all possible ways of completing the $+/\cdot$ operation tables of $0, 1, a$ in such a way that all of the field axioms are satisfied.  It~is a fun exercise (with a lot of similarities to playing Sudoku) to show that the only possible configurations are those given in the left side of Figure~\ref{f:ffoptables}.  In~fact, one can easily check that these tables match those for the well-known field $\ZZ_3$.  
In general, this approach works for any prime value~$p$ to produce the finite field $\ZZ_p$.  


When constructing finite fields of non prime cardinality, such as $4$ (a prime power), we can use the same approach.  Let us consider elements $\{0, 1, a, b\}$.  After checking all the ways to fill out the $+$/$\cdot$ tables, we see there is again a unique solution, depicted in Figure~\ref{f:ffoptables}.  As there is only one way to complete the tables, we once again obtain a unique finite field of this size, which we denote by~$\FF_4$.  Note that $\FF_4$ has characteristic~$2$ (i.e., $1+1 = 0$), so $\FF_4$ is \textbf{not} simply $\ZZ_4$ (which is, in particular, not a field).  

\begin{figure}
\begin{center}
\begin{tabular}{c}
\text{ } \\ 
\begin{tabular}{@{\,\,}c@{\,\,}|@{\,\,}c@{\,\,}|@{\,\,}c@{\,\,}|@{\,\,}c@{\,\,}|}
$+$ & $0$ & $1$ & $a$ \\
\hline
$0$ & $0$ & $1$ & $a$ \\
\hline
$1$ & $1$ & $a$ & $0$ \\
\hline
$a$ & $a$ & $0$ & $1$ \\
\hline
\end{tabular}
\hspace{0.5cm}
\begin{tabular}{@{\,\,}c@{\,\,}|@{\,\,}c@{\,\,}|@{\,\,}c@{\,\,}|@{\,\,}c@{\,\,}|}
$\cdot$ & $0$ & $1$ & $a$ \\
\hline
$0$ & $0$ & $0$ & $0$ \\
\hline
$1$ & $0$ & $1$ & $a$ \\
\hline
$a$ & $0$ & $a$ & $1$ \\
\hline
\end{tabular}
\end{tabular}
\hspace{0.5cm}
\begin{tabular}{@{\,\,}c@{\,\,}|@{\,\,}c@{\,\,}|@{\,\,}c@{\,\,}|@{\,\,}c@{\,\,}|@{\,\,}c@{\,\,}|}
$+$ & $0$ & $1$ & $a$ & $b$ \\
\hline
$0$ & $0$ & $1$ & $a$ & $b$ \\
\hline
$1$ & $1$ & $0$ & $b$ & $a$ \\
\hline
$a$ & $a$ & $b$ & $0$ & $1$ \\
\hline
$b$ & $b$ & $a$ & $1$ & $0$ \\
\hline
\end{tabular}
\hspace{0.5cm}
\begin{tabular}{@{\,\,}c@{\,\,}|@{\,\,}c@{\,\,}|@{\,\,}c@{\,\,}|@{\,\,}c@{\,\,}|@{\,\,}c@{\,\,}|}
$\cdot$ & $0$ & $1$ & $a$ & $b$ \\
\hline
$0$ & $0$ & $0$ & $0$ & $0$ \\
\hline
$1$ & $0$ & $1$ & $a$ & $b$ \\
\hline
$a$ & $0$ & $a$ & $b$ & $1$ \\
\hline
$b$ & $0$ & $b$ & $1$ & $a$ \\
\hline
\end{tabular}
\end{center}
\caption{Unique operation tables for $\FF_3$ (left) and $\FF_4$ (right).}
\label{f:ffoptables}
\end{figure}

It would be nice to identify $\FF_4$ as a more ``familiar'' ring, as we did with $\FF_3 \cong \ZZ_3$.  One way to do this is to view the elements $0, 1, a, b \in \FF_4$ as the elements $\ol 0,\ol 1,\ol z,\ol{z+1}$ in the quotient ring $\ZZ_2[z]/\langle z^2 + z + 1 \rangle$.  In particular, $\ol{z^2+z+1} = \ol 0$ in this quotient ring, meaning $\ol{z^2} = \ol{-z-1}$.  As such, each element can be represented by a polynomial in $z$ with coefficients in $\ZZ_2$, and terms of degree $2$ and higher can be eliminated via the substitution $\ol{z^2} = \ol{z+1}$, e.g., 
$$(\ol{z + 1})(\ol{z + 1}) = \ol{z^2 + 2z + 1} = \ol{(z + 1) + 2z + 1} = \ol{3z + 2} = \ol z,$$
or equivalently using division by $z^2 + z + 1$, e.g., 
$$(\ol{z + 1})(\ol{z + 1}) = \ol{z^2 + 2z + 1} = \ol 1(\ol{z^2 + z + 1}) + \ol z = \ol z.$$
Using this collection of elements, the $+/\cdot$ tables are then obtained by performing polynomial addition and multiplication, reducing coefficients modulo $2$, and then performing polynomial long division by $z^2 + z + 1$; see Figure~\ref{f:quotientoptables}.   Note that in order for this quotient ring to be a field, it is imperative that $z^2 + z + 1$ is irreducible. (In fact, it is the only degree-2 irreducible polynomial in $\ZZ_2[z]$.)  To verify $z^2 + z + 1$ is indeed irreducible, note that any reducible degree 2 polynomial has a degree-1 factor, and therefore has a root by the Root Theorem.  However, neither element of $\ZZ_2$ is a root of $z^2 + z + 1$, so it must be irreducible.  

\begin{figure}
\begin{center}
\begin{tabular}{c|c|c|c|c|}
+ & $\ol 0$ & $\ol 1$ & $\ol z$ & $\ol{z+1}$ \\
\hline
$\ol 0$ & $\ol 0$ & $\ol 1$ & $\ol z$ & $\ol{z+1}$ \rule{0pt}{2.5ex} \\
\hline
$\ol 1$ & $\ol 1$ & $\ol 0$ & $\ol{z+1}$ & $\ol z$ \rule{0pt}{2.5ex} \\
\hline
$\ol z$ & $\ol z$ & $\ol{z+1}$ & $\ol 0$ & $\ol 1$ \rule{0pt}{2.5ex} \\
\hline  
$\ol{z+1}$ & $\ol{z+1}$ & $\ol z$ & $\ol 1$ & $\ol 0$ \rule{0pt}{2.5ex} \\
\hline
\end{tabular}
\hspace{1cm}
\begin{tabular}{c|c|c|c|c|}
$\cdot$ & $\ol 0$ & $\ol 1$ & $\ol z$ & $\ol{z+1}$ \\
\hline
$\ol 0$ & $\ol 0$ & $\ol 0$ & $\ol 0$ & $\ol 0$ \rule{0pt}{2.5ex} \\
\hline
$\ol 1$ & $\ol 0$ & $\ol 1$ & $\ol z$ & $\ol{z+1}$ \rule{0pt}{2.5ex} \\
\hline
$\ol z$ & $\ol 0$ & $\ol z$ & $\ol{z+1}$ & $\ol 1$ \rule{0pt}{2.5ex} \\
\hline
$\ol{z+1}$ & $\ol 0$ & $\ol{z+1}$ & $\ol 1$ & $\ol z$ \rule{0pt}{2.5ex} \\
\hline
\end{tabular}
\end{center}
\caption{Operation tables for $\ZZ_2[z]/\langle z^2 + z + 1 \rangle$.}
\label{f:quotientoptables}
\end{figure}

As a final example, we construct the finite field of $8$ elements, $\FF_8$ (a similar illustration of the construction of $\FF_9$ can be found in~\cite{biggs}).  
Proceeding as above, we wish to use polynomial quotient rings to write $\FF_8$ in the form $\ZZ_2[z]/\langle f(z) \rangle$, where $f$ is an irreducible polynomial of degree $3$.  By the Root Theorem, any reducible degree-3 polynomial has a root, so by inspection of all $2^3 = 8$ polynomials of degree $3$ in $\ZZ_2[z]$, we see that only two are irreducible, namely $z^3 + z + 1$ and $z^3 + z^2 + 1$.  Let
$$\FF_8 = \ZZ_2[z]/\langle z^3 + z + 1 \rangle \qquad \text{and} \qquad \FF_8' = \ZZ_2[w]/\langle w^3 + w^2 + 1 \rangle.$$
Although the two quotient rings above are both fields with $8$ elements, their multiplication ``rules'' appear different, in that in $\FF_8$ we reduce terms of degree $3$ and higher using the equality $\ol{z^3 + z + 1} = \ol 0$, while in $\FF_8'$ we reduce using $\ol{w^3 + w^2 + 1} = \ol 0$.  
For~example, despite the visual similarity, the left-hand side products 
\begin{align}
\label{eq:F8ops1}
\FF_8 &: (\ol{z^2 + 1})(\ol{z + 1}) \bmod (\ol{z^3 + z + 1}) = \ol{z^2}, \\
\label{eq:F8ops2}
\FF_8' &: (\ol{w^2 + 1})(\ol{w + 1}) \bmod (\ol{w^3 + w^2 + 1}) = \ol w
\end{align}
yield visually distinct results.  That said, $\FF_8$ and $\FF_8'$ are isomorphic by the FTFF, and we give an explicit isomorphism in Section~\ref{sec:mainproof}.  

\section{The Key Lemma:\ factoring over finite fields}
\label{sec:keylemma}

For a finite field $\FF_q$ with $q = p^r$, the key to the FTFF turns out to be factoring the polynomial $x^q - x$, both over $\FF_q$ itself and over $\ZZ_p$.  The Key Lemma, stated below and followed immediately by several examples, identifies precisely how $x^q - x$ factors as a product of irreducible polynomials over both fields.  

\begin{keylemma}
Suppose $q = p^r$ for $p$ prime and $r \in \ZZ_{\ge 1}$.  
\begin{enumerate}[(a)]
\item 
If $\KK$ is any finite field with $|\KK| = q$, then the polynomial $x^q - x$ factors over $\KK$ as a product of distinct linear factors.  

\item 
The polynomial $x^q - x$ factors over $\ZZ_p$ as the product of all irreducible polynomials over $\ZZ_p$ with degree dividing $r$.  

\end{enumerate}
\end{keylemma}

Let us work through a few examples.  We start with $q = 4$, which we constructed as 
$$\FF_4 = \ZZ_2[z]/\langle z^2 + z + 1 \rangle = \{\ol 0, \ol 1, \ol z, \ol{z + 1}\}$$
in Section~\ref{sec:constructing}.  Since $0, 1 \in \ZZ_2$ are both roots of $x^4 - x$, the Root Theorem tells us $x$ and $x - 1$ are both factors.  The Key Lemma implies all remaining factors are degree~2.  Since polynomial long division by $x - 1$ yields
$$x^4 - x = x(x - 1)(x^2 + x + 1),$$
the Key Lemma implies $x^2 + x + 1$ is the \textbf{only} irreducible polynomial of degree~$2$ over~$\ZZ_2$, a fact we also observed in Section~\ref{sec:constructing}.  
Now, since $\ZZ_2 \subsetneq \FF_4$, the ``extra'' two elements of $\FF_4$ provide more coefficients at our disposal when factoring, so some irreducible polynomials over $\ZZ_2$, like $x^2 + x + 1$ in this case, may be factored further over $\FF_4$.   Indeed, we obtain four distinct linear factors, one for each element of $\FF_4$, i.e.,
\begin{align*}
x^4 - x
&= x(x - 1)(x^2 + x + 1)
\\
&= x(x - \ol 1)(x - \ol z)(x - \ol{z+1})
\end{align*}
wherein $\ol z$ and $\ol{z+1}$ are the roots of $x^2 + x + 1$ in $\FF_4$.  
Remember that the polynomials in the above expression live in $\FF_4[x]$, so in the second line $\ol z$ and $\ol{z+1}$ are \textbf{coefficients} that live in~$\FF_4$.  

For $q = 8$, after factoring $x$ and $x + 1$ out of $x^8 - x$, we obtain a degree-6 polynomial that, by the Key Lemma, must factor into (exactly 2) distinct degree-3 irreducible factors over $\ZZ_3$.  As both irreducible polynomials were identified in Section~\ref{sec:constructing}, this~yields
\begin{align*}
x^8 - x
&= x(x + 1)(x^6 + x^5 + x^4 + x^3 + x^2 + x + 1)
\\
&= x(x + 1)(x^3 + x + 1)(x^3 + x^2 + 1)
\end{align*}
as the factorization over $\ZZ_3$.   It is in this way that we can use the Key Lemma to locate \textbf{all} possible choices of an irreducible polynomial when constructing a finite field of a particular size.  
Next, to factor further over $\FF_8$, we choose the representation $\FF_8 = \ZZ_2[z]/\langle z^3 + z + 1 \rangle$ and obtain
\begin{align*}
x^3 + x + 1
&= (x - \ol{z})(x - \ol{z^2})(x - \ol{z^2+z})
\\
x^3 + x^2 + 1
&= (x - \ol{z+1})(x - \ol{z^2+1})(x - \ol{z^2+z+1}).
\end{align*}
Had we instead chosen the representation $\FF_8' = \ZZ_2[w]/\langle w^3 + w^2 + 1 \rangle$, we would have obtained
\begin{align*}
x^3 + x + 1
&= (x - \ol{w+1})(x - \ol{w^2+1})(x - \ol{w^2+w})
\\
x^3 + x^2 + 1
&= (x - \ol{w})(x - \ol{w^2})(x - \ol{w^2+w+1})
\end{align*}
as the remaining linear factors of $x^8 - x$.  

We give one final example before proving the Key Lemma.   Applying a similar process as above for $q = 9$, we obtain
$$x^9 - x = x(x + 1)(x + 2)(x^2 + 1)(x^2 + x + 2)(x^2 + 2x + 2)$$
over $\ZZ_3$, and factoring further over $\FF_9 = \ZZ_3[z]/\langle z^2 + 1 \rangle$
(this time there were 3 possible representations to choose from)
yields
\begin{align*}
x^2 + 1
&= (x - \ol z)(x - \ol{2z})
\\
x^2 + x + 2
&= (x - \ol{z+1})(x - \ol{2z+1})
\\
x^2 + 2x + 2
&= (x - \ol{z+2})(x - \ol{2z+2}),
\end{align*}
which we encourage the reader to verify as an exercise.  


\begin{proof}[Proof of the Key Lemma]
Suppose $q = p^r$ for $p$ prime and $r \ge 1$, and suppose $\KK$ is a field with $|\KK| = q$.  
Since $\KK$ is a field, $(\KK \setminus \{0\}, \cdot)$ is a group of order $q - 1$, meaning that every element has order dividing $q-1$.  As~such, $a^{q-1} - 1 = 0$ for every nonzero $a \in \KK$, and thus each is a root of $x^q - x$.  By the Root Theorem, this produces $q$ distinct linear factors and $x^q - x$ has degree $q$, so this must be precisely the list of factors, proving part~(a).  

Next, fix an irreducible polynomial $f \in \ZZ_p[x]$, and let $d = \deg f$.  We wish to show $f(x) \mid x^q - x$ if and only if $d \mid r$, as this implies that the irreducible factors of $f(x)$ over $\ZZ_p$ claimed in part~(b) are precisley those that appear.  Since $f$ is irreducible, as stated in Section \ref{sec:background}, the quotient ring $K = \ZZ_p[x]/\langle f(x) \rangle$ is a field.  To more clearly distinguish $K$ from the field $\FF_q$ constructed elsewhere in this document, we will denote the elements of $K$ using the ``bracket'' notation $[h(x)]$ for $h \in \ZZ_p[x]$ rather than with the ``overline'' notation.  Since $|K| = p^d$, we can list the elements of $K$ as
$$K = \big\{ [h_1(x)], \ldots, [h_{p^d}(x)] \big\}$$
with $h_1(x) = 0$.  If $f$ is linear, then the claim follows from part~1, so assume $d \ge 2$.  

First, suppose $d \mid r$.  Since $K$ is a field, multiplication is cancellative, so multiplying by $[x]$ permutes the set of nonzero elements.  In particular, the list
$$[x][h_2(x)], \quad [x][h_3(x)], \quad \ldots, \quad [x][h_{p^d}(x)]$$
contains every nonzero element of $K$ exactly once.  As such, the product of all elements in this list can be simplified in two ways to obtain
\begin{align*}
[x][h_2(x)] \cdots [x][h_{p^d}(x)]
&= [h_2(x)][h_3(x)] \cdots [h_{p^d}(x)]
\\
&= [x^{p^d-1}][h_2(x)] \cdots [h_{p^d}(x)],
\end{align*}
in $K$, 
where the expressions on either side of the first equality consist of the product of the same $p^d - 1$ elements of $K$ (albeit in a different order).  
Subtracting and factoring yields
$$[x^{p^d-1} - 1][h_2(x)][h_3(x)] \cdots [h_{p^d}(x)] = [0] \in K,$$
which implies $[x^{p^d-1} - 1] = [0]$ since $K$ has no zero-divisors.  This means we have $x^{p^d-1} - 1 \in \langle f(x) \rangle \subseteq \ZZ_p[x]$ and thus $f(x) \mid x^{p^d-1} - 1$.  Since $d \mid r$, say $r = dk$ for some $k \in \ZZ$, 
$$p^r - 1 = (p^d - 1)(p^{d(k-1)} + p^{d(k-2)} + \cdots + p^d + 1),$$
meaning $p^d - 1 \mid p^r - 1$.  Analogously, fixing $t \in \ZZ_{\ge 0}$ so that $p^r - 1 = t(p^d - 1)$, we have
$$x^{q - 1} - 1 = x^{t(p^d - 1)} - 1 = (x^{p^d - 1} - 1)(x^{(p^d-1)(t-1)} + x^{(p^d-1)(t-2)} + \cdots + 1).$$
Putting all of this together, we conclude $f(x) \mid x^{p^d-1} - 1 \mid x^q - x$.  

Conversely, suppose $f(x) \mid x^q - x$.  Using the division algorithm to write $r = ad + b$ for $a, b \in \ZZ$ with $0 \le b < d$, we wish to show $b = 0$.  By way of contradiction, suppose that $b$ is positive.  Since $|K| = p^d$, similar reasoning as in the first paragraph of this proof implies $[x]^{p^d} = [x]$ in $K$, and by assumption $[x^q - x] = [0]$ in $K$, so 
$$[x] = [x^q] = [x^{p^{ad + b}}] = [((\cdots((\underbrace{x^{p^d})^{p^d})\cdots)^{p^d}}_{a \text{ times}})^{p^b}] = [x^{p^b}].$$
By the Freshman's Dream, for any $g(x) = g_0 + g_1x + g_2x^2 + \cdots \in \ZZ_p[x]$, we have
\begin{align*}
[g(x)]^{p^b}
&= [(g_0)^{p^b} + (g_1)^{p^b}(x^{p^b}) + (g_2)^{p^b}(x^{p^b})^2 + \cdots] \\
&= [g_0 + g_1x + g_2x^2 + \cdots] \\
&= [g(x)],
\end{align*}
meaning every element of $K$ is a root of $x^{p^b} - x$.  However, this is impossible by the Root Theorem since $K$ has $p^d > p^b$ elements, so we conclude $b = 0$.  This completes the proof that $f(x) \mid x^q - x$ if and only if $d \mid r$.  

There remains one final claim to prove:\ that each irreducible polynomial $f(x)$ in the factorization of $x^q - x$ over $\ZZ_p$ appears only once.  Indeed, by part~(a), the roots of $x^q - x$ in $K$ are all distinct, so $x^q - x$ cannot have repeated factors over $\ZZ_p$ as this would yield repeated roots in $K$.  This completes the proof.  
\end{proof}

\section{The fundamental theorem}
\label{sec:mainproof}

In this section, we use the Key Lemma to prove the FTFF in its entirety.  Before diving into the proof, let's briefly explore some of its implications in the context of $\FF_8$ and $\FF_8'$ from Section~\ref{sec:constructing}.  

First, the set $\FF_8 \setminus \{\ol 0\}$ is ensured to be a cyclic group under multiplication, meaning there is some element $a \in \FF_8$ such that the list $a, a^2, a^3, \ldots$ includes every nonzero element of $\FF_8$.  One such element is $\ol{z + 1}$, and we can readily check that every nonzero element of $\FF_8$ can be written as $(\ol{z + 1})^n$ for some $n$.  In fact, it turns out that any nonzero element we choose for $a$ will do the trick (except $a = \ol 1$, of course).  This is not true in general:\ in $\FF_7$ (which is isomorphic to $\ZZ_7$), only 2 nonzero elements generate $\{1,2,3,4,5,6\}$ as a group under multiplication modulo $7$ (one such element is $3 \in \ZZ_7$, and we encourage the reader to locate the other).  

Second, the FTFF implies that $\FF_8$ and $\FF_8'$ are isomorphic, but it is not hard to show that the map $\FF_8 \rightarrow \FF_8'$ given by $\ol{a z^2 + b z + c} \mapsto \ol{a w^2 + b w + c}$ is \textbf{not} an isomorphism (compare, for instance, the right hand sides of~\eqref{eq:F8ops1} and~\eqref{eq:F8ops2} in Section~\ref{sec:constructing}).  One possible isomorphism $\sigma: \FF_8 \rightarrow \FF^{'}_8$ turns out to be
$$\begin{array}{r@{}c@{}l@{\qquad}r@{}c@{}l@{\qquad}r@{}c@{}l}
\sigma(\ol 0) &{}={}& \ol 0 &
\sigma(\ol z) &{}={}& \ol{w + 1} &
\sigma(\ol{z^2 + z}) &{}={}& \ol{w^2 + w} \\
\sigma(\ol 1) &{}={}& \ol 1 &
\sigma(\ol{z^2}) &{}={}& \ol{w^2 + 1} &
\sigma(\ol{z^2 + z + 1}) &{}={}& \ol{w^2 + w + 1} \\
\sigma(\ol{z + 1}) &{}={}& \ol w &
\sigma(\ol{z^2 + 1}) &{}={}& \ol{w^2},
\end{array}$$
which happens to map a generator to another generator. In~general, locating an explicit isomorphism between finite fields of equal size need not be easy, as mapping a generator to a generator does not always yield an isomorphism.  This map sends the element $\ol{z + 1} \in \FF_8$ to the element $\ol{w} \in \FF_8'$, and the remaining nonzero elements, necessarily of the form $(\ol{z + 1})^n$ for some $n \ge 2$ by the previous paragraph, is sent to $(\ol{w})^n$.  This~guarantees multiplication is preserved by $\sigma$.  Verifying that addition is preserved can be done manually, or by observing that every element $a \in \FF_8$ can be written uniquely as a sum involving $\ol 1$, $\ol z$, and $\ol{z^2}$, and that for each such $a$, $\sigma(a)$ equals precisely the image of this sum (for example, $\sigma(\ol{z^2 + 1}) = \sigma(\ol{z^2}) + \sigma(\ol 1) = \ol{w^2}$).  


\begin{proof}[Proof of the FTFF]
Consider the subring $R \subset \FF_q$ consisting of $0, 1, 1 + 1, \ldots \in \FF_q$, and let $p = |R|$ (the \emph{characteristic} of $\FF_q$).  We see that $R \iso \ZZ_p$, and so $p$ must be prime, as otherwise $\ZZ_p$ (and thus $\FF_q$) would contain zero-divisors.  This makes~$\FF_q$ a vector space over the field $\ZZ_p$, necessarily finite dimensional since $\FF_q$ is finite, so the fundamental theorem of linear algebra tells us that for some $r \ge 1$, we have $(\FF_q, +) \iso (\ZZ_p)^r$.  This proves part~(c)(i) and the forward direction of part~(a).  

For the backwards direction of part~(a), we must prove $\ZZ_p[x]/\<f(x)\>$ is a field with exactly $p^r$ elements. Thus, it is enough to show there exists at least one irreducible polynomial in $\ZZ_p[x]$ of each degree $r$.  The Key Lemma implies that the sum of the degrees of all irreducible polynomials in $\ZZ_p[x]$ whose degree divides $r$ is $p^r$.  If we sum only those degrees strictly dividing~$r$, we obtain
$$\sum_{d \mid r, \, d \ne r} p^d \le \sum_{d < r} p^d = \frac{p^r - 1}{p - 1} < p^r.$$
As such, there is an irreducible polynomial $f \in \ZZ_p[x]$ with $\deg f = r$, as~desired. 

Next, we prove part~(c)(ii).  Let $N$ denote the maximum order of any element of the group $(\FF_q \setminus \{0\}, \cdot)$.  We claim every element of $(\FF_q \setminus \{0\}, \cdot)$ has order dividing $N$.  Indeed,~if $|a| = N$ and $|b| = m \nmid N$, then there exists some prime power $t$ such that $t \mid m$ and $t \nmid N$.  However, $|ab^{m/t}| = \lcm(N,t) > N$ contradicts the maximality of $N$.  This proves the claim.  Now, this means every nonzero element of $\FF_q$ is a root of $x^N - 1$, which is only possible if $\deg(x^N - 1) \ge q - 1 = |\FF_q \setminus \{0\}|$.  As such, $N = q - 1$, and any element of order $N$ generates $(\FF_q \setminus \{0\}, \cdot)$, thereby proving part~(c)(ii).  

Finally, we prove part~(b).  Fix any irreducible polynomial $f \in \ZZ_p[x]$ of degree~$r$.  We claim $\FF_q \iso \ZZ_p[x]/\<f(x)\>$.  Since $f(x)$ divides $x^q - x$ by The Key Lemma, some element $a \in \FF_q$ is a root of $f$.  Consider the homomorphism
$$\begin{array}{r@{}c@{}l}
\varphi:\ZZ_p[x] &{}\longrightarrow{}& \FF_q \\
g(x) &{}\longmapsto{}& g(a),
\end{array}$$
which has kernel 
$$\ker(\varphi) = \{g(x) : g(a) = 0\} = \langle f(x) \rangle$$
by the Root Theorem since $f$ is irreducible over $\ZZ_p$ and has $a$ as a root.  As such, the First Isomorphism Theorem implies $\ZZ_p[x]/\langle f(x) \rangle \cong \image(\varphi)$, and $\varphi$ must be surjective since $\FF_q$ and $\image(\varphi)$ both have $q$~elements, so the claimed isomorphism is shown.  
\end{proof}

\bibliography{finitefields}
\bibliographystyle{amsplain}

\end{document}